\documentclass[12pt]{amsart}
\usepackage{amsmath, amssymb, amsthm, bm, float, mathtools, enumitem, caption}
\usepackage{hyperref}

\usepackage[noabbrev,capitalize]{cleveref}
\usepackage{adjustbox}
\crefname{equation}{}{}
\usepackage{ytableau}
\usepackage{graphicx, tikz}
\usepackage[textheight=8.5in, textwidth=6.75in]{geometry}
\usepackage{subcaption}
\usepackage{microtype}
\usepackage{stmaryrd}
\usepackage{mathtools}
\newtheorem{theorem}{Theorem}[section]
\newtheorem{lemma}[theorem]{Lemma}

\newtheorem{proposition}[theorem]{Proposition}

\newtheorem*{conjecture*}{Conjecture}

\theoremstyle{definition}

\theoremstyle{remark}

\numberwithin{equation}{section}

\title[Parity biases for missing and present integers in partitions]{Parity biases for missing and present integers in partitions}
\date{\today}
\author{Julia Q.D. Du \and Gurinder Singh}
\address{School of Mathematical Sciences, Hebei Key Laboratory of Computational Mathematics and Applications, Hebei Research Center of the Basic
discipline Pure Mathematics, Hebei Normal University, Shijiazhuang 050024, P.R. China}
\email{qddu@hebtu.edu.cn}
\address{Postdoctoral Research Station of Mathematics, Hebei Normal University, Shijiazhuang 050024, P.R. China}
\email{gurindermaan1018@gmail.com}

\begin{document}

\begin{abstract}
Bhoria, Eyyunni and Santra proved that, for each positive integer $k$, the total number of distinct positive integers not exceeding the largest part and occurring fewer than $k$ times in partitions of $n$ equals the total number of parts different from $k$. We establish an overpartition analogue in which the corresponding total is twice a statistic defined by the parity of the nonoverlined parts, together with an overlined-part contribution when $k$ is even. For $k=1$, we compare the parity distributions of missing integers and present parts. We prove an alternating parity bias for the number of parts greater than one in ordinary partitions, determine the resulting strict ordering of the four parity classes for $n>34$, and establish two congruences modulo $2$. For overpartitions, the difference between the numbers with an even and an odd number of nonoverlined even parts equals the number of overpartitions into odd parts. Combining explicit coefficient estimates with previously established positivity results, we obtain a strict ordering of the four overpartition parity classes for $n>27$.
\end{abstract}

\maketitle
%%%%%%%%%%%%%%%%%%%%%%%%%%%%%%%%%%%%%%%%%%%%%%%%%%_Introduction_%%%%%%%%%%%%%%%%%%%%%%%%%%%%%%%%%%%%%%%%%%%%%%%%%%%%%%%%%%%%%%%%%%%%%%%%%%%%%%%%%%%%%%%%%%%%%%%%%%%%%%%%%%%%%%%%%%%%%%%%%%%%%%%%%%%%%%%%%%%%%%%%%%%%%%%%%%%%%%%%%%%%%%%%%%%%%%%%%%%%%%%%%%%%%%%%%%%%%%%%%%%%%%%%%%%%%%%%%%%%%%%%%%%%%%%%%%%%%
\section{Introduction}
Mathematicians have been interested in partition statistics such as the parts that occur and the integers that are absent in partitions. 
Statistics associated with the parts occurring in a partition have been studied extensively. Andrews \cite{And_2008} introduced the smallest-parts function, $\operatorname{spt}(n)$, which sums the multiplicities of the smallest parts over all partitions of $n$, and established congruences for this function. Another direction concerns the total number of parts in restricted partition families. Proving conjectures of Beck, Andrews \cite{And_2017} showed that the excess of the total number of parts in partitions of $n$ into odd parts over the corresponding total for partitions into distinct parts equals the number of partitions of $n$ containing exactly one even part size, which may appear with arbitrary positive multiplicity.
Ballantine and Bielak \cite{BB_2019} subsequently gave combinatorial proofs of these identities. Parity questions also arise naturally: Lehmer’s identity expresses the difference between the numbers of partitions having an even and an odd number of even parts as the number of partitions into distinct odd parts. Ballantine et al. \cite{BBF_2022} established generalizations and companion identities for this result.
The missing-integers viewpoint includes the \emph{minimal excludant}, or mex, of a partition: the smallest positive integer that is not a part. Grabner and Knopfmacher \cite{GK_2006} studied this statistic under the name ``least gap.'' Andrews and Newman \cite{AN_2019} subsequently developed its combinatorial theory and studied $\sigma\mathrm{mex}(n)$, the sum of the minimal excludants over all partitions of $n$. A combinatorial interpretation was given by Ballantine and Merca \cite{BM_2021}. Connections with the crank were investigated by Hopkins, Sellers and Stanton \cite{HSS_2022}. Baruah et al. \cite{BEM_2023} refined the arithmetic function $\sigma\mathrm{mex}(n)$ by considering the sums of odd and even minimal excludants separately. Some related works are \cite{BS_2021, CR_2021, SS_2020}. At the other end of the set of missing integers, Chern \cite{Ch_2021} introduced the \emph{maximal excludant}, the largest nonnegative integer below the largest part that does not occur. Bhoria, Eyyunni and Santra \cite{BES_2026} studied the entire set of missing positive integers below the largest part, rather than only its extremal elements. They obtained generating functions, congruences, and parity-bias conjectures for partitions and overpartitions. Aniceto and Ballantine \cite{AB_2026} later gave combinatorial proofs of several of their identities. The authors \cite{DS_2026} proved the three conjectures posed in \cite{BES_2026}.
\par We use the notation
\begin{align*}
(a;q)_r=\prod_{j=0}^{r-1}(1-aq^j),\qquad
 (a;q)_\infty=\prod_{j=0}^{\infty}(1-aq^j),\qquad f_r(q)=(q^r;q^r)_\infty.
\end{align*}
Infinite products and analytic identities are considered for $|q|<1$; coefficient identities also hold as formal power series. An overpartition is a partition in which the first occurrence of each distinct part size may be overlined. Let $p(n)$ and $\overline{p}(n)$ denote the ordinary partition and overpartition functions, respectively. In particular,
\begin{equation*}
 \overline{P}(q):=\sum_{n=0}^{\infty}\overline p(n)q^n
       =\prod_{j=1}^{\infty}\frac{1+q^j}{1-q^j}=\frac{(-q;q)_{\infty}}{(q;q)_{\infty}}.
\end{equation*}
%%%%%%%%%%%%%%%%%%%%%%%%%%%%
The following result of Bhoria et al. \cite{BES_2026} provides the starting point for our work.
\begin{theorem}\label{Thm_BES} \cite[Theorem 2.9]{BES_2026}
For a positive integer $k$, let $\mathcal{M}(k,n)$ be the number of distinct positive integers not exceeding the largest part and appearing less than $k$ times in all partitions of $n$, and let $\mathcal{P}(k,n)$ be the number of parts different from $k$ in all partitions of $n$. Then, for all nonnegative integers $n$,
\begin{align*}
\mathcal{M}(k,n)=\mathcal{P}(k,n).
\end{align*}
\end{theorem}
We answer the question in \cite{BES_2026} asking for an overpartition analogue of Theorem \ref{Thm_BES}. For a positive integer $k$, let $\overline{\mathcal{M}}(k,n)$ be the number of distinct positive integers not exceeding the largest part and appearing less than $k$ times in all overpartitions of $n$. Set
\begin{align}\label{eqSSk}
\mathcal{S}_k:=
\{1\leq j<k:j\equiv k \pmod2\}\cup\{j>k:j\not\equiv k\pmod2\}.   
\end{align}
Let $\mu_k(n)$ be the number of nonoverlined parts whose sizes belong to $\mathcal{S}_k$ in all overpartitions of $n$, and $\overline{\mu}(n)$ be the number of overlined parts in all overpartitions of $n$. Let
\begin{align*}
    \overline{\mathcal{P}}(k,n)=
    \begin{cases}
        \mu_k(n), &\text{if } k \text{ is odd};\\
        \mu_k(n)+\overline{\mu}(n), & \text{if } k \text{ is even}.
    \end{cases}
\end{align*}
%We obtain the following an overpartition analogue of Theorem \ref{Thm_BES}.

\begin{theorem}\label{Thm_main_1}
 For every $k\geq1$ and $n\geq0$,
\begin{align*}
\overline{\mathcal{M}}(k,n)=2\;\overline{\mathcal{P}}(k,n).
\end{align*}
In particular, the total number of missing positive integers below the largest part in overpartitions of $n$ is twice the total number of nonoverlined even parts.
\end{theorem}
%%%%%%%%%%%%%%%%%%%%%%%%%%%%%%%%%%%%%%%
Let $\mathcal{M}_e(n)$ (resp., $\mathcal{M}_o(n)$) denote the number of partitions of $n$ with an even (resp., odd) number of missing integers less than the largest part. Let $\mathcal{P}_e(n)$ (resp., $\mathcal{P}_o(n)$) denote the number of partitions of $n$ in which the parts greater than $1$ are even (resp., odd) in number. Then,
\begin{align}\label{M_P_relation}
\mathcal{M}_e(n)+\mathcal{M}_o(n)=p(n)=\mathcal{P}_e(n)+\mathcal{P}_o(n).
\end{align}
Bhoria, Eyyunni and Santra \cite{BES_2026} proved
\begin{align}\label{M_gf}
\sum_{n=0}^{\infty}\left(\mathcal{M}_e(n)-\mathcal{M}_o(n)\right)q^n=\frac{(-2q;q)_{\infty}}{(-q;q)_{\infty}}.
\end{align}
They conjectured and the authors \cite{DS_2026} proved that for $n>34$, 
\begin{align}\label{eq1.1}
\mathcal{M}_e(n) >\mathcal{M}_o(n).
\end{align}
%%%%%%%%%%%%%%%%%%%%%%%%%%%%%%%%%%%%%%%%%%%%%%%%%%%
For $k=1$, Theorem \ref{Thm_BES} proves that the number of missing integers in all the partitions of $n$
equals the number of parts different from 1 in all the partitions of $n$. Although these statistics have the same total over all partitions of $n$, this equality does not determine their respective parity
distributions. It is therefore natural to investigate whether there is a comparable bias in the parity of the number of parts different from 1 in all partitions. We obtain the following alternating bias for the present-part statistic beyond the same threshold.
\begin{theorem}\label{Thm_main_2}
For even integers $n>34$, we have
$$\mathcal{P}_e(n)>\mathcal{P}_o(n).$$
For odd integers $n>34$, we have
$$\mathcal{P}_e(n)<\mathcal{P}_o(n).$$
\end{theorem}
Theorem \ref{Thm_main_2} and \eqref{eq1.1} motivate us to compare $\mathcal{M}_e(n)$, $\mathcal{M}_o(n)$, $\mathcal{P}_e(n)$, and $\mathcal{P}_o(n)$ on a single scale. The comparison gives us the following.
\begin{theorem}\label{Thm_main_3}
For even $n>34$, we have
\begin{align}\label{MP-ineq-even}
\mathcal{M}_e(n)>\mathcal{P}_e(n)>\mathcal{P}_o(n)>\mathcal{M}_o(n),
\end{align}
and for odd $n>34$, we have
\begin{align}\label{MP-ineq-odd}
\mathcal{M}_e(n)>\mathcal{P}_o(n)>\mathcal{P}_e(n)>\mathcal{M}_o(n).
\end{align}
\end{theorem}
Moreover, we prove that $\mathcal{M}_e(n)$ and $\mathcal{P}_e(n)$ are equal modulo $2$, and so are $\mathcal{M}_o(n)$ and $\mathcal{P}_o(n)$.
\begin{theorem}\label{Thm_main_3a}
For every $n\geq0$, we have
\begin{align*}
\mathcal{M}_e(n)\equiv \mathcal{P}_e(n) \pmod{2},\qquad\mathcal{M}_o(n)\equiv \mathcal{P}_o(n) \pmod{2}.
\end{align*}
\end{theorem}
%%%%%%%%%%%%%%%%%%%%%%%%%%%%%%%%%%%%%%%%%%%%%%%%%%
Now, we enter the world of overpartitions. Let $\overline{\mathcal{M}}_e(n)$ (resp., $\overline{\mathcal{M}}_o(n)$) denote the number of overpartitions of $n$ with an even (resp., odd) number of missing integers less than the largest part. It is natural to be curious about the corresponding partition functions $\overline{\mathcal{P}}_e(n)$ and $\overline{\mathcal{P}}_o(n)$. Motivated by Theorem \ref{Thm_main_1} with $k=1$, we define $\overline{\mathcal{P}}_e(n)$ (resp., $\overline{\mathcal{P}}_o(n)$) to be the number of overpartitions of $n$ with an even (resp., odd) number of nonoverlined even parts. Then,
\begin{align}\label{eqmfe}
\overline{\mathcal{M}}_e(n)+\overline{\mathcal{M}}_o(n)=\overline{p}(n)=\overline{\mathcal{P}}_e(n)+\overline{\mathcal{P}}_o(n).
\end{align}
The generating function proved in \cite{BES_2026} is
\begin{align*}
\sum_{n=0}^{\infty}\left(\overline{\mathcal{M}}_e(n)-\overline{\mathcal{M}}_o(n)\right)q^n=\frac{(-3q;q)_{\infty}}{(-q;q)_{\infty}}.
\end{align*}
The authors of \cite{BES_2026} also conjectured the following overpartition analogue of \eqref{eq1.1}: For $n>27$,
\begin{align}\label{eq_oM}
\overline{\mathcal{M}}_e(n) >\overline{\mathcal{M}}_o(n).
\end{align}
We also proved this conjecture in \cite{DS_2026}.  
We show that the present-part parity difference in the case of overpartitions has a direct positive interpretation.
\begin{theorem}\label{Thm_main_4}
Let $\overline{\mathcal{P}}_{\text{odd}}(n)$ denote the number of overpartitions of $n$ into odd parts. Then for all positive integers $n$, we have
\begin{align*}
\overline{\mathcal{P}}_e(n) -\overline{\mathcal{P}}_o(n)=\overline{\mathcal{P}}_{\text{odd}}(n)>0.
\end{align*}
\end{theorem}
The next natural step is to compare $\overline{\mathcal{M}}_e(n)$, $\overline{\mathcal{M}}_o(n)$, $\overline{\mathcal{P}}_e(n)$, and $\overline{\mathcal{P}}_o(n)$ as in Theorem \ref{Thm_main_3}. 
\begin{theorem}\label{Thm_main_5}
For $n>27$, we have
\begin{align*}
\overline{\mathcal{P}}_e(n)>\overline{\mathcal{M}}_e(n)>\overline{\mathcal{M}}_o(n)>\overline{\mathcal{P}}_o(n).
\end{align*}
\end{theorem}
%%%%%%%%%%%%%%%%%%%%%%%%
The rest of this paper is organized as follows. Section \ref{Sec2} proves Theorem \ref{Thm_main_1}. Section \ref{Sec3} establishes Theorem \ref{Thm_main_2} by a coefficient argument based on Sylvester's identity. Section \ref{Sec4} proves Theorems \ref{Thm_main_3} and \ref{Thm_main_3a}; it explicitly states the quantitative estimate used from \cite{DS_2026}. Section \ref{Sec5} proves Theorems \ref{Thm_main_4} and \ref{Thm_main_5}. We conclude with questions for $k\ge2$.
%%%%%%%%%%%%%%%%%%%%%%%%%%%%%%%%%%%%%%%%%%%%%%%%%%%%%%%%%%%%%%%%%%%%%%%%%%%%%%%%%%%%%%%%%%%%%%%%%%%%
\section{Proof of Theorem \ref{Thm_main_1}}\label{Sec2}
For a partition $\lambda$, let $\ell(\lambda)$ and $a(\lambda)$ denote the number of parts (length of $\lambda$) and the largest part of $\lambda$, respectively. We denote the total number of occurrences of size $j$ (multiplicity of $j$) in $\lambda$ by $m_j(\lambda)$. The multiplicity of a part of an overpartition includes both its overlined and nonoverlined occurrences. For the empty partition, set $\ell(\lambda)=a(\lambda)=0$. Let $\overline{\mathcal{P}}(n)$ be the set of all overpartitions of $n$. For a property $\mathfrak{p}$, $\mathbf{1}_{\mathfrak{p}}$ is $1$ if $\mathfrak{p}$ is true and $0$ otherwise.

\begin{proof}[Proof of Theorem \ref{Thm_main_1}]
% From the definition of $\overline{\mathcal{M}}(k,n)$, we have
% \begin{align}\label{eq2.1}
% \overline{\mathcal{M}}(k,n)=\sum_{\lambda\in\overline{\mathcal P}(n)}
% \#\{1\le j\le a(\lambda):m_j(\lambda)<k\},    
% \end{align}
% where integers with multiplicity zero are included. 

Since for any $\lambda\in\overline{\mathcal P}(n)$,
\begin{align*}
\#\{1\le j\le a(\lambda):m_j(\lambda)<k\}
=
a(\lambda)-\sum_{j=1}^{\infty}
\mathbf{1}_{\{m_j(\lambda)\ge k\}},  
\end{align*}
we have
\begin{align*}
\overline{\mathcal{M}}(k,n)=
\sum_{\lambda\in\overline{\mathcal P}(n)} \#\{1\le j\le a(\lambda):m_j(\lambda)<k\}
= \sum_{\lambda\in\overline{\mathcal P}(n)}\ell(\lambda)
-\sum_{\lambda\in\overline{\mathcal P}(n)}\sum_{j=1}^{\infty}
\mathbf{1}_{\{m_j(\lambda)\ge k\}},
\end{align*}
where we used the fact that, by conjugation, the sum of the largest parts equals the total number of parts.
%here we used the total of the largest parts equals the total number of parts by conjugation.
Then
\begin{align}\label{eq2.2}
\sum_{n=0}^{\infty}\overline{\mathcal{M}}(k,n)q^n=\sum_{n=0}^{\infty}L(n)q^n
-\sum_{n=0}^{\infty}O_k(n)q^n,
\end{align}
where
\begin{align*}
L(n):=\sum_{\lambda\in\overline{\mathcal P}(n)}\ell(\lambda),\qquad O_k(n):=\sum_{\lambda\in\overline{\mathcal P}(n)}\sum_{j=1}^{\infty}
\mathbf{1}_{\{m_j(\lambda)\ge k\}}.
\end{align*}
\par Let $L(n,m)$ count the number of overpartitions of $n$ with length $m$. Then
$L(n)=\sum_{m=0}^{\infty}mL(n,m)$.
Since the bivariate generating function for the sequence $L(n,m)$ is given by
\begin{align*}
\mathcal{L}(z;q):=\sum_{n=0}^{\infty}\sum_{m=0}^{\infty} L(n,m)z^mq^n=\prod_{r=1}^{\infty}\frac{1+zq^r}{1-zq^r},
\end{align*}
we have
\begin{align}\label{eq2.3}
\sum_{n=0}^{\infty}L(n)q^n=\frac{\partial}{\partial z}\left(\mathcal{L}(z;q)\right)\Bigg\vert_{z=1}=2\overline{P}(q)\sum_{j=1}^{\infty}\frac{q^j}{1-q^{2j}}.
\end{align}
\par Let $O_k(n,m)$ count the number of overpartitions of $n$ having exactly $m$ different part sizes with multiplicity at least $k$. Then $O_k(n)=\sum_{m=0}^{\infty}mO_k(n,m)$.
Notice that the bivariate generating function for $O_k(n,m)$ is given by
\begin{align*}
\mathcal{O}_k(z;q):=\sum_{n=0}^{\infty}\sum_{m=0}^{\infty} O_k(n,m)z^mq^n=\prod_{r=1}^{\infty}\left(1+2\sum_{s=1}^{k-1}q^{rs}+2z\sum_{s=k}^{\infty}q^{rs}\right),
\end{align*}
so we have
\begin{align}\label{eq2.4}
\sum_{n=0}^{\infty}O_k(n)q^n=\frac{\partial}{\partial z}\left(\mathcal{O}_k(z;q)\right)\Bigg\vert_{z=1}=2\overline{P}(q)\sum_{j=1}^{\infty}\frac{q^{kj}}{1+q^{j}}.
\end{align}
Combining \eqref{eq2.2}, \eqref{eq2.3}, and \eqref{eq2.4}, we obtain
\begin{align}\label{eq2.5}
\sum_{n=0}^{\infty}\overline{\mathcal{M}}(k,n)q^n=2\overline{P}(q)\sum_{j=1}^{\infty}\left(\frac{q^j}{1-q^{2j}}-\frac{q^{kj}}{1+q^{j}}\right).
\end{align}
It is easy to see that
\begin{align*}
\frac{q^j}{1-q^{2j}}-\frac{q^{kj}}{1+q^{j}}=\sum_{r\in\mathcal{S}_k}q^{jr}+\mathbf{1}_{2\mid k}\frac{q^j}{1+q^j},
\end{align*}
where $\mathcal{S}_k$ is defined in \eqref{eqSSk}.
Thus, \eqref{eq2.5} becomes
\begin{align}\label{eq2.6}
\sum_{n=0}^{\infty}\overline{\mathcal{M}}(k,n)q^n=2\overline{P}(q)\left(
\sum_{r\in\mathcal{S}_k}\frac{q^r}{1-q^r}+\mathbf{1}_{2\mid k}\sum_{j=1}^{\infty}\frac{q^j}{1+q^j}\right).
\end{align}
%%%%%%%%%%%%%%%%%%%%%%%%%%%%%%%%%%%
\par For a positive integer $r$, let $\nu_r(n)$ denote the number of nonoverlined parts equal to $r$ in all overpartitions of $n$. Since
\begin{align*}
\overline{\mathcal{P}}(k,n)=\mu_k(n)+\mathbf{1}_{2\mid k}\overline{\mu}(n),
\end{align*}
we have
\begin{align}\label{eq2.7}
\overline{\mathcal{P}}(k,n)%=\mu_k(n)+\mathbf{1}_{2\mid k}\overline{\mu}(n)
=\sum_{r\in\mathcal{S}_k}\nu_r(n)+\mathbf{1}_{2\mid k}\overline{\mu}(n).
\end{align}
Let $V_r(n,m)$ be the number of overpartitions of $n$ with $m$ nonoverlined parts equal to $r$. Then $\nu_r(n)=\sum_{m=0}^{\infty}mV_r(n,m)$. 
Since
\begin{align*}
\mathcal{V}_r(z;q):=\sum_{n=0}^{\infty}\sum_{m=0}^{\infty} V_r(n,m)z^mq^n=
\frac{1+q^r}{1-zq^r}\cdot\prod_{\substack{j=1\\ j\neq r}}^{\infty}\frac{1+q^j}{1-q^j},
\end{align*}
we have
\begin{align}\label{eq2.8}
\sum_{n=0}^{\infty}\nu_r(n)q^n=\frac{\partial}{\partial z}\left(\mathcal{V}_r(z;q)\right)\Bigg\vert_{z=1}=\overline{P}(q)\frac{q^{r}}{1-q^{r}}.
\end{align}
Let $\overline{U}(n,m)$ be the number of overpartitions of $n$ with $m$ overlined parts. Then, 
\begin{align*}
\overline{\mathcal{U}}(z;q):=\sum_{n=0}^{\infty}\sum_{m=0}^{\infty} \overline{U}(n,m)z^mq^n=
\prod_{j=1}^{\infty}\frac{1+zq^j}{1-q^j}.
\end{align*}
Similarly, $\overline{\mu}(n)=\sum_{m=0}^{\infty}m\overline{U}(n,m)$ implies that
\begin{align}\label{eq2.9}
\sum_{n=0}^{\infty}\overline{\mu}(n)q^n=\frac{\partial}{\partial z}\left(\overline{\mathcal{U}}(z;q)\right)\Bigg\vert_{z=1}=\overline{P}(q)\sum_{j=1}^{\infty}\frac{q^{j}}{1+q^{j}}.
\end{align}
From \eqref{eq2.7}, \eqref{eq2.8}, and \eqref{eq2.9}, we obtain
\begin{align}\label{eq2.10}
\sum_{n=0}^{\infty}\overline{\mathcal{P}}(k,n)q^n=\overline{P}(q)\left(
\sum_{r\in\mathcal{S}_k}\frac{q^r}{1-q^r}+\mathbf{1}_{2\mid k}\sum_{j=1}^{\infty}\frac{q^j}{1+q^j}\right).
\end{align}
We conclude from \eqref{eq2.6} and \eqref{eq2.10} that
\begin{align*}
\overline{\mathcal{M}}(k,n)=2\;\overline{\mathcal{P}}(k,n).  
\end{align*}
This completes the proof of Theorem \ref{Thm_main_1}.
\end{proof}

%%%%%%%%%%%%%%%%%%%%%%%%%%%%%%%%%%%%%%%%%%%%%%%%%%%%%%%%%%%%%%%%%%%%%%%%%%%%%%%%%%%%%%%%%%%%%%%%%%%%%%%%%%%%%%%%%%%%%%%%%%%%%%%%%%%%%%%%%%%%%%%%%%%%Proof_of_Theorem_2_%%%%%%%%%%%%%%%%%%%%%%%%%%%%%%%%%%%%%%%%%%%%%%%%%%%%%%%%%%%%%%%%%%%%%%%%%%%%%%%%%%%%%%%%%%%%%%%%%%%%%%%%%%%%%%%%%%%%%%%%%%%%%%%%%%%%%%%%%%%%%%%%%%%%%%%%%%%%%%%%%%%%%%%%%%%%%%%%%%%%%%%%%%%%%%%%%%%%%%%%%%%%%%%%%%%%%%%%%%%%%%%%%%%%%%%%%%%%%%%%%%%%%%%%%%%%
\section{Proof of Theorem \ref{Thm_main_2}}\label{Sec3}
Let $b(m,n)$ be the number of partitions of $n$ with exactly $m$ parts greater than $1$. Then, 
\begin{align*}
\sum_{n=0}^{\infty}\sum_{m=0}^{\infty}b(m,n)z^m q^n=\frac{1}{1-q}\cdot\frac{1}{(zq^2;q)_{\infty}}.
\end{align*}
Taking $z=-1$, we have
\begin{align}\label{P_gf}
G(q):=\sum_{n=0}^\infty (\mathcal{P}_e(n)-\mathcal{P}_o(n))q^n
&=\sum_{n=0}^{\infty}\sum_{m=0}^{\infty}b(m,n)(-1)^m q^n\nonumber\\
&=\frac{1}{1-q}\cdot\frac{1}{(-q^2;q)_{\infty}}
=(1+q)(q^3;q^2)_\infty,
\end{align}
where we used Euler's product identity $1/(-q;q)_\infty=(q;q^2)_\infty$ 
in the last equality.
 
Let
\begin{align}\label{eqw2}
\sum_{n=0}^{\infty}g(n)q^n:=G(-q)=(1-q)(-q^3;q^2)_{\infty}.
\end{align}
Then
$$g(n)=(-1)^n\left(\mathcal{P}_e(n)-\mathcal{P}_o(n)\right).$$
Thus, Theorem \ref{Thm_main_2} is equivalent to 
the following statement.
%$g(n)>0$, for all $n>34$. That is,
\begin{theorem}\label{Thm_7.4}
For any $n>34$, we have 
\begin{align*}
g(n)>0.
\end{align*}
\end{theorem}
For $\alpha(q)\in\mathbb{Z}[[q]]$, we use $[q^m]\alpha(q)$ to denote the coefficient of $q^m$ in the series $\alpha(q)$.
Recall an identity of Sylvester \cite[Eq. (9.2.3)]{And_1976}:
\begin{align}\label{eqw1}
(-xq;q)_{\infty}=\sum_{k=0}^{\infty}\frac{(-xq;q)_k}{(q;q)_k}(1+xq^{2k+1})x^kq^{k(3k+1)/2}.
\end{align}
Replacing $x$ with $1/q$ and $q$ with $q^2$ in \eqref{eqw1}, and then dividing both sides by $1+q$, we get
\begin{align}\label{eqw3}
(-q^3;q^2)_{\infty}=1+\sum_{k=1}^{\infty}\frac{q^{3k^2}(1+q^{4k+1})(-q^3;q^2)_{k-1}}{(q^2;q^2)_k}.
\end{align}
Define, for $k\geq1$,
\begin{align*}
T_k(q):=\frac{(1-q)(-q^3;q^2)_{k-1}}{(q^2;q^2)_k},\qquad
Q_k(q):=(1+q^{4k+1})T_k(q),\qquad R(q):=\frac{1}{1-q^2}-q.
\end{align*}
Then, multiplying \eqref{eqw3} by $1-q$ and using \eqref{eqw2}, we get
\begin{align}\label{eqw1a}
\sum_{n=0}^\infty g(n)q^n =1-q+\sum_{k=1}^{\infty}q^{3k^2}Q_k(q).
\end{align}

% Multiplying both sides of \eqref{eqw3} with $(1-q)$ yields
% \begin{align}\label{eqw1a}
% G(-q)=1-q+\sum_{k=1}^{\infty}q^{3k^2}Q_k(q).
% \end{align}
To prove the positivity of the coefficients $g(n)$, we need the following three lemmas for the coefficients of $T_k(q)$ and $Q_k(q)$.
%%%%%%%%%%%%%%%%%%%%%%%%%%%%%%%%%%%%_Lemma_1_%%%%%%%%%%%%%%%%%%%%%%%%%%%%%%
\begin{lemma}\label{sign-P-R}
For any $k\geq3$ and $n\geq 0$,
\begin{align}\label{eqw8}
[q^n]\left(T_k(q)-R(q)\right)\geq 0.
\end{align}
\end{lemma}

\begin{proof}
We prove \eqref{eqw8} by induction on $k$.
For $k=3$, we have
\begin{align*}
T_3(q)=\frac{(1-q)(1+q^3)(1+q^5)}{(1-q^2)(1-q^4)(1-q^6)}
=\frac{(1-q+q^2)(1+q^5)}{(1-q^4)(1-q^6)}.
\end{align*}
Subtracting $R(q)$ from $T_3(q)$, we obtain
\begin{align}\label{eqw3a}
T_3(q)-R(q)=\frac{q^8+q^{11}}{(1-q^4)(1-q^6)},
\end{align}
which implies that $T_3(q)-R(q)$ has nonnegative coefficients.

Assume that $T_k(q)-R(q)$ has nonnegative coefficients
for a fixed integer $k\geq 3$.
From the definition of $T_k(q)$, we obtain
\begin{align}\label{eqw4}
T_{k+1}(q)=\frac{1+q^{2k+1}}{1-q^{2k+2}}T_k(q).
\end{align}
Subtracting $R(q)$ from both sides of \eqref{eqw4} gives
\begin{align}\label{eqw5}
T_{k+1}(q)-R(q)=\frac{(1+q^{2k+1})(T_k(q)-R(q))+q^{2k+1}(1+q)R(q)}{1-q^{2k+2}}.
\end{align}
Notice that
\begin{align*}
(1+q)R(q)=\frac{1}{1-q}-q-q^2=1+\frac{q^3}{1-q},
\end{align*}
which, together with the induction hypothesis and \eqref{eqw5}, implies that $T_{k+1}(q)-R(q)$ has nonnegative coefficients. Hence, by induction it proves \eqref{eqw8} for all $k\geq3$.
\end{proof}
%%%%%%%%%%%%%%%%%%%%%%%%%%%%%%%%%%%%_Lemma_2_%%%%%%%%%%%%%%%%%%%%%%%%%%%%%%
\begin{lemma}\label{lemma_7.2}
For any $k\geq3$ and $n\geq 0$, 
$[q^n](Q_k(q)+q)\geq 0$.
\end{lemma}

\begin{proof}
Since
\begin{align*}
Q_k(q)+q&=(Q_k(q)-(1+q^{4k+1})R(q))+(1+q^{4k+1})R(q)+q\\
&=(1+q^{4k+1})(T_k(q)-R(q))+(1+q^{4k+1})R(q)+q\\
&=(1+q^{4k+1})(T_k(q)-R(q))+\left(\frac{1}{1-q^2}-q^{4k+2}\right)+\frac{q^{4k+1}}{1-q^2},
\end{align*}
by Lemma \ref{sign-P-R}, we conclude that
$Q_k(q)+q$ has nonnegative coefficients.
\end{proof}
%%%%%%%%%%%%%%%%%%%%%%%%%%%%%%%%%%%%_Lemma_3_%%%%%%%%%%%%%%%%%%%%%%%%%%%%%%
\begin{lemma}\label{lemma_7.3}
For every integer $m\geq15$, $[q^m]Q_3(q)\geq2$.
\end{lemma}
\begin{proof}
Multiplying both sides of \eqref{eqw3a} by $1+q^{13}$, we get
\begin{align}\label{eqw11}
Q_3(q)=(1+q^{13})R(q)+\frac{(1+q^{13})(q^8+q^{11})}{(1-q^4)(1-q^6)}:=K_1(q)+K_2(q),
\end{align}
where
$$K_1(q)=(1+q^{13})R(q)=\frac{1}{1-q^2}+\frac{q^{13}}{1-q^2}-q-q^{14}$$
and
$$K_2(q)=\frac{(1+q^{13})(q^8+q^{11})}{(1-q^4)(1-q^6)}.$$
It is clear that for $m\geq 15$,
\begin{align}\label{eqw12}
[q^m]K_1(q)=1.
\end{align}
For $K_2(q)$, observe that
\begin{align*}
\frac{1}{(1-q^4)(1-q^6)}=\sum_{a,b\geq0}q^{4a+6b}.
\end{align*}
Notice that 
every even integer $n\geq4$ can be written as $4i+6j$ with $i,j\geq0$ (if $4\mid n$, take $i=n/4,~j=0$; 
otherwise take $i=(n-6)/4,~j=1$). 
So, for even $m\geq16$,
$$[q^m]\left(\frac{q^8}{(1-q^4)(1-q^6)}\right)\geq1,$$
and for odd $m\geq15$,
$$[q^m]\left(\frac{q^{11}}{(1-q^4)(1-q^6)}\right)\geq1.$$
Therefore, for $m\geq 15$,
\begin{align}\label{eqw13}
[q^m]K_2(q)\geq 1.
\end{align}
Combining \eqref{eqw11}, \eqref{eqw12}, and \eqref{eqw13}, we conclude that for all $m\geq15$,
$[q^m]Q_3(q)\geq2$.
The proof is complete.
\end{proof}
%%%%%%%%%%%%%%%%%%%%%%%%%%%%%%%%%%%%_Theorem_%%%%%%%%%%%%%%%%%%%%%%%%%%%%%%
Now, we are ready to prove the positivity of the coefficients $g(n)$.

\begin{proof}[Proof of Theorem \ref{Thm_7.4}]

%Recall that from \eqref{eqw1a}
%\begin{align*}
%\sum_{n=0}^\infty g(n)q^n =1-q+\sum_{k=1}^{\infty}q^{3k^2}Q_k(q).
%\end{align*}
We explicitly discuss the sign of the expression
for the values of $k$ appearing in the sum
on the R.H.S. of \eqref{eqw1a}.
For $k=1$,
\begin{align*}
q^3Q_1(q)=q^3\frac{1+q^5}{1+q}=q^3-q^4+q^5-q^6+q^7,
\end{align*}
which implies that for $n\geq 8$,
\begin{align}\label{eqw14}
[q^n](1-q+q^3Q_1(q))=0.
\end{align}

For $k=2$,
\begin{align*}
q^{12}Q_2(q)&=q^{12}\frac{(1-q)(1+q^3)(1+q^9)}{(1-q^2)(1-q^4)}\\
% &=q^{12}\frac{(1-q+q^2)(1+q^9)}{1-q^4}\\
% &=q^{12}\left(\frac{-q+q^2+q^9-q^{10}}{1-q^4}+\frac{1+q^{11}}{1-q^4}\right)\\
% &=q^{12}\left(-q+q^2-q^5+q^6+\frac{1+q^{11}}{1-q^4}\right)\\
&=-q^{13}+q^{14}-q^{17}+q^{18}+\frac{q^{12}+q^{23}}{1-q^4},
\end{align*}
so, for $n\geq 18$, 
\begin{align}\label{eqw15}
[q^n](q^{12}Q_2(q))\geq0.
\end{align}

For $k=3$, it follows from Lemma \ref{lemma_7.3} that
for $n\geq 42$,
\begin{align}\label{eqw16}
[q^n](q^{27}Q_3(q))=[q^{n-27}]Q_3(q)\geq2.
\end{align}

Lemma \ref{lemma_7.2} implies that for any $n\geq 0$,
\begin{align}\label{eqw17}
[q^n]\left(\sum_{k=4}^{\infty}q^{3k^2}Q_k(q)+\sum_{k=4}^{\infty}q^{3k^2+1}\right)
\geq 0.
\end{align}
Since the exponents $3k^2+1$ ($k\geq4$) are pairwise distinct, every coefficient of $\sum_{k\geq4}q^{3k^2+1}$ is at most 1. Thus \eqref{eqw16} and \eqref{eqw17} imply that
for $n\geq 42$,
\begin{align}\label{eqw18}
[q^n]\left(q^{27}Q_3(q)+\sum_{k=4}^{\infty}q^{3k^2}Q_k(q)\right)\geq1.
\end{align}
Combining \eqref{eqw1a}, \eqref{eqw14}, \eqref{eqw15}, and \eqref{eqw18}, we conclude that for $n\geq 42$,
\begin{align*}
[q^n](G(-q))=[q^n]\left(1-q+\sum_{k=1}^{\infty}q^{3k^2}Q_k(q)\right)\geq1.
\end{align*}
This implies that for $n\geq 42$, $g(n)>0$. 
\par It remains to check the seven values of $n$, precisely $35\leq n\leq 41$. One can easily check that $g(35)=3$, $g(36)=1$, $g(37)=1$, $g(38)=1$, $g(39)=3$, $g(40)=2$, and $g(41)=1$. 
Hence, the proof is complete.
\end{proof}

%%%%%%%%%%%%%%%%%%%%%%%%%%%%%%%%%%%%%%%%%%%%%%%%%%%%%%%%%%%%%%%%%%%%%%%%%%%%%%%%%%%%%%%%%%%%%%%%%%%%%%%%%%%%%%%%%%%%%%%%%%%%%%%%%%%%%%%%%%%%%%%%%%%%%%%%%%%%%%%%%%%%%%%%%%%%%%%%%%%%%%%%%%%%%%%%%%%%%%%%%%%%%%%%%%%%%%%%%%%%%%%%%%%%%%%%%%%%%%%%%%%%%%%%%%%%%%%%%%%%%%%%%%%%%%%%%%%%%%%%%%%%%%%%%%%%%%%%%%%%%%%%%%%%%%%%%%%%%%%%%%%%%%%%%%%%%%%%%%%%%%%%%%%%%%%%%%%%%%%%%%%%%%%%%%%%%%%%%%%%%%%%%%%%%%%%%%%%%%%%%%%%%%%%%%%%%%%%%%%%%%%%%%%%%%%%%%%%%%%%%%%%%%%%%%%%%%%%%%%%%%%%%%%%%%%%%%%%%%%%%%%%%%%%%%%%%%%%%%%%%%%%%%%%%%%%%%%%%%%%%%%%%%%%%%%%%%%%%%%%%%%%%%%%%%%%%%%%%%%%%%%%%%%%%%%%%%%%%%%%%%%%%%%%%%%%%%%%%%%%%%%%%%%%%%%%%%%%%%%%%%%%%%%%%%%%%%%%%%%%%%%%%%%%%%%%%%%%%%%%%%%%%%%%%%%%%%%%%%%%%%%%%%%%%%%%%%%%%%%%%%%%%%%%%%%%%%%%%%%%%%%%%%%%%%%%%%%%%%%%%%%%%%%%%%%%%%%%%%%%%%%%%%%%%%%%%%%%%%%%%%%%%%%%%%%%%%%%
%%%%%%%%%%%%%%%%%%%%%%%%%%%%%%%%%%%%%%%%%%%%%%%%%%%%%%%%%%%%%%%%%%%%%%%%%%%%%%%%%%%%%%%%%%%%%%%%%%%%%%%%%%%%%%%%%%%%%%%%%%%%%%%%%%%%%%%%%%%%%%%%%%%%Proof_of_Theorem_3_%%%%%%%%%%%%%%%%%%%%%%%%%%%%%%%%%%%%%%%%%%%%%%%%%%%%%%%%%%%%%%%%%%%%%%%%%%%%%%%%%%%%%%%%%%%%%%%%%%%%%%%%%%%%%%%%%%%%%%%%%%%%%%%%%%%%%%%%%%%%%%%%%%%%%%%%%%%%%%%%%%%%%%%%%%%%%%%%%%%%%%%%%%%%%%%%%%%%%%%%%%%%%%%%%%%%%%%%%%%%%%%%%%%%%%%%%%%%%%%%%%%%%%%%%%%%
\section{Proof of Theorems \ref{Thm_main_3} and \ref{Thm_main_3a}}\label{Sec4}
%%%%%%%%
\par Let 
\begin{align}\label{eq4.5}
a(n):=\mathcal{M}_e(n)-\mathcal{M}_o(n)   
\end{align}
and
\begin{align}\label{eq4.6}
C(q)=\sum_{n=0}^{\infty}c(n)q^n:=(-q^3;q^2)_{\infty}=\prod_{j=1}^{\infty}(1+q^{2j+1}).
\end{align}
We first prove a coefficient comparison stronger than is needed for Theorem \ref{Thm_main_3}. 
We use the following explicit estimate from the proof of the coefficient-positivity theorem in \cite[Theorem 1.4]{DS_2026}: For $n\ge8846$, 
one has
\begin{align}\label{eq5.1}
a(n)>\frac{1-e^{-1}}{2\pi\sqrt {6A_2}} x^{3/2}e^{2A_2/x},
\end{align}
where
\begin{align*}
 A_2=\frac{(\log2)^2}{2}+\frac{\pi^2}{12}+\operatorname{Li}_2(-1/2),
 \qquad x=\sqrt{A_2/n}.
\end{align*}
For $|z|<1$, the dilogarithm is defined by  
\begin{align*}
\operatorname{Li}_2(z)=\sum_{m=1}^{\infty}\frac{z^m}{m^2}.
\end{align*}
Its principal branch extends analytically to $\mathbb C\setminus[1,\infty)$ and has the integral representation
\begin{align*}
\text{Li}_2(z)=-\int_{0}^{z}\frac{\operatorname{Log}(1-u)}{u}\,du.
\end{align*}
Here, $\operatorname{Log}$ denotes the principal complex logarithm. See \cite{Zagier} or \cite[Section 2]{DS_2026} for more information on this function. The alternating series gives
$$-\frac{65}{144}=-\frac12+\frac1{16}-\frac1{72}<\operatorname{Li}_2(-1/2)<-\frac{1}{2}+\frac1{16}=-\frac{7}{16}.$$
This, along with the elementary bounds $2/3<\log2<7/10$ and $49/5<\pi^2<10$, yields
\begin{align}\label{eqA_2bd}
\frac9{16}<A_2<1,\qquad
\frac{1-e^{-1}}{2\pi\sqrt6}\,A_2^{1/4}>\frac1{64}.    
\end{align}
%For the second inequality, we use $1-e^{-1}>1/2$, $A_2^{1/4}>4/5$, and $2\pi\sqrt6<16$.
%%%%%%%%%%%%%%%
\begin{lemma}\label{lemma4}
For all $n>34$, we have
\begin{align*}
a(n)>c(n).
\end{align*}
\end{lemma}

\begin{proof}
Substituting $x=\sqrt{A_2/n}$ into \eqref{eq5.1} and using \eqref{eqA_2bd}, we obtain
\begin{align}\label{eq5.2}
a(n)>\frac1{64}n^{-3/4}e^{\frac32\sqrt n},\qquad n\ge8846.
\end{align}
For an upper bound for $c(n)$, let $h(s)=\log(1+e^{-s})$.
%Recall that
%\begin{align*}
%C(q):=\sum_{n\ge0}c(n)q^n
%=\prod_{j\ge1}(1+q^{2j+1}).    
%\end{align*}
%Let $g(s)=\log(1+e^{-s})$.
Since
\begin{align*}
h''(s)=\frac{e^{-s}}{(1+e^{-s})^2}>0,   
\end{align*}
the midpoint inequality on $[2jt,(2j+2)t]$ gives
\begin{align*}
2t\,h((2j+1)t)
\le \int_{2jt}^{(2j+2)t}h(s)\,ds.  
\end{align*}
Consequently, for $t>0$,
\begin{align*}
\log C(e^{-t})\le\sum_{j=0}^{\infty}h((2j+1)t)\le\frac1{2t}\int_0^\infty\log(1+e^{-s})\,ds=\frac{1}{2t}\sum_{r=1}^\infty \frac{(-1)^{r+1}}{r^2}=\frac{\pi^2}{24t},
\end{align*}
where we used
\begin{align*}
\sum_{r=1}^\infty \frac{(-1)^{r+1}}{r^2}=\frac{\pi^2}{12}.
\end{align*}
Since $c(n)\ge0$,
\begin{align*}
c(n)e^{-nt}\le C(e^{-t})
\le e^{\pi^2/(24t)}.    
\end{align*}
Choosing $t=\pi/\sqrt{24n}$, we find
\begin{align}\label{eq5.3}
c(n)\le e^{\pi\sqrt{n/6}}
<e^{\frac43\sqrt n}\qquad(n\ge1).   
\end{align}
%%%%%%%%%%%%%%%
\par To compare \eqref{eq5.2} and \eqref{eq5.3}, define
\begin{align*}
E(u)=\frac{\sqrt u}{6}-\frac34\log u-\log64.    
\end{align*}
Then $E'(u)=(\sqrt u-9)/(12u)>0$ for $u>81$. Moreover,
\begin{align*}
E(8100)=15-\frac34\log8100-\log64>15-\frac{30}{4}-5=\frac52>0.
\end{align*}
Therefore, $E(n)>0$ for $n\ge8846$, which means
\begin{align*}
\frac1{64}n^{-3/4}e^{\frac32\sqrt n}
>e^{\frac43\sqrt n}.
\end{align*}
Combining this with \eqref{eq5.2} and \eqref{eq5.3} proves that for $n\ge8846$, $a(n)>c(n)$.
For $35\le n\le8845$, we verify $a(n)>c(n)$ using Mathematica. This completes the proof. 
\end{proof}
%%%%%%%%%%%%%%%%%%%%%%%%%%%%%%%%%%%%%%%%%%%%%%%%%%%%%%%%%%%%%%%%
\begin{proof}[Proof of Theorem \ref{Thm_main_3}]
From \eqref{P_gf} and \eqref{eq4.6}, we have for all $n\geq1$, 
\begin{align}\label{eq4.7}
\mathcal{P}_e(n)-\mathcal{P}_o(n)=(-1)^n(c(n)-c(n-1)),
\end{align}
which, together with \eqref{M_P_relation} and \eqref{eq4.5}, gives
\begin{align*}
\frac{a(n)-c(n)+c(n-1)}{2}=\begin{cases}
\mathcal{M}_e(n)-\mathcal{P}_e(n), &n\text{ even},\\
\mathcal{M}_e(n)-\mathcal{P}_o(n),&n\text{ odd}.
\end{cases}
\end{align*}
It follows from Lemma \ref{lemma4} and $c(n)\geq 0$ that
for all even $n>34$, $\mathcal{M}_e(n)>\mathcal{P}_e(n)$,
and for all odd $n>34$, $\mathcal{M}_e(n)>\mathcal{P}_o(n)$.
This, together with \eqref{M_P_relation} and Theorem \ref{Thm_main_2},
gives \eqref{MP-ineq-even} and \eqref{MP-ineq-odd}.
\end{proof}
%%%%%%%%%%%%%%%%%%%%%%%%%%%%%%%%%%%%%%%%%%%%%%%%%%%%%%%%%%%%%%%%
\begin{proof}[Proof of Theorem \ref{Thm_main_3a}]
From \eqref{M_P_relation}, we obtain
\begin{align}\label{eq4.1}
2(\mathcal{M}_e(n)-\mathcal{P}_e(n))=\mathcal{M}_e(n)-\mathcal{M}_o(n)-\mathcal{P}_e(n)+\mathcal{P}_o(n).
\end{align}
Subtracting \eqref{P_gf} from \eqref{M_gf} and using \eqref{eq4.1} yields
\begin{align}\label{eq4.2}
2\sum_{n=0}^{\infty}\left(\mathcal{M}_e(n)-\mathcal{P}_e(n)\right)q^n=\frac{1}{(-q;q)_{\infty}}\left((-2q;q)_{\infty}-\frac{1+q}{1-q}\right),
\end{align}
where we used Euler's theorem $(-q;q)_{\infty}=1/(q;q^2)_{\infty}$. Since
\begin{align*}
(-2q;q)_{\infty}=\prod_{j=1}^{\infty}(1+2q^j)\equiv 1+2\sum_{j=1}^{\infty}q^j=\frac{1+q}{1-q}\pmod{4},
\end{align*}
\eqref{eq4.2} implies that
\begin{align}\label{eq4.3}
\mathcal{M}_e(n)\equiv \mathcal{P}_e(n)\pmod2.   
\end{align}
Combining \eqref{M_P_relation} and \eqref{eq4.3}, we get that
\begin{align*}
\mathcal{M}_o(n)\equiv \mathcal{P}_o(n)\pmod2.   
\end{align*}
This completes the proof of Theorem \ref{Thm_main_3a}.
\end{proof}
%%%%%%%%%%%%%%%%%%%%%%%%%%%%%%%%%%%%%%%%%%%%%%%%%%%%%%%%%%%%%%%%%%%%%%%%%%%%%%%%%%%%%%%%%%%%%%%%%%%%%%%%%%%%%%%%%%%%%%%%%%%%%%%%%%%%%%%%%%%%%%%%%%%%Proof_of_Theorem_1.6_%%%%%%%%%%%%%%%%%%%%%%%%%%%%%%%%%%%%%%%%%%%%%%%%%%%%%%%%%%%%%%%%%%%%%%%%%%%%%%%%%%%%%%%%%%%%%%%%%%%%%%%%%%%%%%%%%%%%%%%%%%%%%%%%%%%%%%%%%%%%%%%%%%%%%%%%%%%%%%%%%%%%%%%%%%%%%%%%%%%%%%%%%%%%%%%%%%%%%%%%%%%%%%%%%%%%%%%%%%%%%%%%%%%%%%%%%%%%%%%%%%%%%%%%%%%%
\section{Proof of Theorems \ref{Thm_main_4} and \ref{Thm_main_5}}\label{Sec5}
We first give a proof of Theorem \ref{Thm_main_4}.
\begin{proof}[Proof of Theorem \ref{Thm_main_4}]
Let $\overline{\mathcal{P}}_{\rm{noe}}(n,m)$ count the number of overpartitions of $n$ with $m$ non-overlined even parts. Then,
\begin{align}\label{eq6.1}
\sum_{n=0}^{\infty}\sum_{m=0}^{\infty} \overline{\mathcal{P}}_{\rm{noe}}(n,m)z^mq^n=\prod_{j=1}^{\infty}\frac{1+q^j}{1-q^{2j-1}}\cdot \frac{1}{1-zq^{2j}}.
\end{align}
Since
\begin{align*}
\overline{\mathcal{P}}_e(n)-\overline{\mathcal{P}}_o(n)=\sum_{m=0}^{\infty}(-1)^m\overline{\mathcal{P}}_{\rm{noe}}(n,m),
\end{align*}
by substituting $z=-1$ into \eqref{eq6.1}, we obtain
\begin{align*}
\sum_{n=0}^{\infty}\left(\overline{\mathcal{P}}_e(n)-\overline{\mathcal{P}}_o(n)\right)q^n=\frac{(-q;q^2)_{\infty}}{(q;q^2)_{\infty}}=\sum_{n=0}^{\infty}\overline{\mathcal{P}}_{\rm{odd}}(n)q^n.
\end{align*}
For $n\geq1$, $\overline{\mathcal{P}}_{\rm{odd}}(n) \geq 1$,
since there exists an overpartition of $n$ consisting of $n$ nonoverlined ones.
Notice that $\overline{\mathcal{P}}_{\rm{odd}}(0)=1$.
Thus, $\overline{\mathcal{P}}_{\rm{odd}}(n)>0$ for all $n\geq0$.
\end{proof}
%%%%%%%%%%%%%%%%%%%%%%%%%%%%%%%%%%%%%%%%%%%%%%%%%%%%%%%%%%%%%%%%%%%%%%%%%%%%%%%%%%%%%%%%%%%%%%%%%%%%%%%%%%%%%%%%%%%%%%%%%%%%%%%%%%%%%%%%%%%%%%%%
From \eqref{eqmfe}, we have
\begin{align}\label{eq7.}
\overline{\mathcal{P}}_e(n)-\overline{\mathcal{M}}_e(n)=\overline{\mathcal{M}}_o(n)-\overline{\mathcal{P}}_o(n).
\end{align}
Therefore,
\begin{align*}
\overline{\mathcal{P}}_e(n)>\overline{\mathcal{M}}_e(n) \quad \iff \quad \overline{\mathcal{M}}_o(n)>\overline{\mathcal{P}}_o(n).
\end{align*}
Thus, in order to prove Theorem \ref{Thm_main_5}, that is, to prove that for $n>27$,
\begin{align*}
\overline{\mathcal{P}}_e(n)>\overline{\mathcal{M}}_e(n)>\overline{\mathcal{M}}_o(n)>\overline{\mathcal{P}}_o(n),
\end{align*}
in view of \eqref{eq_oM}, it is enough to prove the following.
\begin{proposition}\label{prop5.1}
For $n>27$,
\begin{align*}
\overline{\mathcal{P}}_e(n)>\overline{\mathcal{M}}_e(n).
\end{align*}
\end{proposition}
Put
\begin{align}\label{eq7.0}
\alpha(n):=\overline{\mathcal{M}}_e(n)-\overline{\mathcal{M}}_o(n),\qquad \beta(n):=\overline{\mathcal{P}}_e(n)-\overline{\mathcal{P}}_o(n),
\end{align}
and
\begin{align}
 \mathcal{A}(q):=\sum_{n=0}^{\infty}\alpha(n)q^n
 &=\frac{(-3q;q)_\infty}{(-q;q)_\infty}=\prod_{j=1}^{\infty}\frac{1+3q^j}{1+q^j},\label{eq7.1}\\
 \mathcal{B}(q):=\sum_{n=0}^{\infty}\beta(n)q^n
 &=\frac{(-q;q^2)_\infty}{(q;q^2)_\infty}=\frac{f_2(q)^3}{f_1(q)^2f_4(q)}.\label{eq7.2}
\end{align}
Combining \eqref{eq7.} and \eqref{eq7.0}, we find that
\begin{align}\label{eq7p}
2\left(\overline{\mathcal{P}}_e(n)-\overline{\mathcal{M}}_e(n)\right)=\beta(n)-\alpha(n).
\end{align}
Thus, in order to prove Proposition \ref{prop5.1},
we only need to prove that $\beta(n)>\alpha(n)$ for $n>27$.
More precisely, we obtain an upper bound for $\alpha(n)$ (in Lemma \ref{lemma5.2}) and a lower bound for $\beta(n)$ (in Lemma \ref{lemma5.3}), and then compare them.
%%%%%%%%%%%%%%%%%%%%%%_Lemma_5.2_%%%%%%%%%%%%%%%%%%%%%%%%%%%%%%%%%%%%%%%%%%%%%%%%
\begin{lemma}\label{lemma5.2}
 For any $n>27$,
\begin{align}\label{eq:a-upper}
\alpha(n)<2 e^{(53/25)\sqrt n}.
\end{align}
\end{lemma}
\begin{proof}
Define
\begin{align}\label{eq7J}
J:=\int_0^\infty \log\left(\frac{1+3e^{-s}}{1+e^{-s}}\right)\,ds .   
\end{align}
The integrand in \eqref{eq7J} is positive and decreasing as its derivative is
$$-\frac{2e^{-s}}{(1+3e^{-s})(1+e^{-s})}<0.$$
Then, substituting $q=e^{-t}$ in \eqref{eq7.1} and
applying the right-endpoint rectangle bound on each interval $[(j-1)t,jt]$, for $t>0$, gives
\begin{align}\label{eq7.4}
\log \mathcal{A}(e^{-t})\le J/t,\qquad\mathcal{A}(e^{-t})\le e^{J/t}.   
\end{align}
The substitution $u=e^{-s}$ in \eqref{eq7J} gives
\begin{align}
J=&\int_{0}^{1}\frac{\log(1+3u)-\log(1+u)}{u}\,du=-\text{Li}_2(-3)+\text{Li}_2(-1)\nonumber\\
&=\frac{(\log3)^2}{2}+\frac{\pi^2}{12}+\operatorname{Li}_2(-1/3),\label{eq7.3}
\end{align}
where we used $\text{Li}_2(-1)=-\pi^2/12$ and the following real inversion identity \cite[p. 8]{Zagier}:
\begin{align*}
\text{Li}_2\left(-1/u\right)+\text{Li}_2(-u)=-\frac{\pi^2}{6}-\frac{1}{2}(\log u)^2\qquad(u>0).
\end{align*}
The alternating series for $\operatorname{Li}_2(-1/3)$ yields
\begin{align*}
-\frac13<\operatorname{Li}_2(-1/3)<-\frac13+\frac1{36}=-\frac{11}{36}.   
\end{align*}
With this, we use $13/12<\log3<11/10$ and $49/5<\pi^2<247/25$ in \eqref{eq7.3} to obtain
\begin{align}\label{eq_J_bd}
1<\frac{1541}{1440}<J<\frac{2021}{1800}
 <\left(\frac{53}{50}\right)^2.
\end{align}
From \eqref{eq_oM}, $\alpha(n)>0$ for $n>27$,
and the only negative coefficients of $\mathcal{A}(q)$ through degree $27$ are
\begin{align*}
\alpha(4)=-2,\quad \alpha(11)=-8,\quad \alpha(13)=-24,\quad \alpha(26)=-96.  
\end{align*}
Therefore, the series
\begin{align*}
\widehat{\mathcal{A}}(q):=\mathcal{A}(q)+2q^4+8q^{11}+24q^{13}+96q^{26}   
\end{align*}
has nonnegative coefficients. For $n\ge28$, its coefficient of $q^n$
is $\alpha(n)$. Consequently, by \eqref{eq7.4},
\begin{align*}
\alpha(n)e^{-nt}\le \widehat{\mathcal{A}}(e^{-t})\le e^{J/t}+130. 
\end{align*}
Choosing $t=\sqrt{J/n}$ yields
\begin{align*}
\alpha(n)\le e^{2\sqrt{Jn}}+130e^{\sqrt{Jn}}.    
\end{align*}
Since $J>1$ and $n\ge28$, we have
$e^{\sqrt{Jn}}>e^{\sqrt{28}}>e^5>130$. From \eqref{eq_J_bd}, $2\sqrt J<53/25$, and therefore
\begin{align*}
\alpha(n)<2e^{2\sqrt{Jn}}<2e^{(53/25)\sqrt n},
\end{align*}
as required. 
\end{proof}
%%%%%%%%%%%%%%%%%%%%%%_Lemma_5.3a_%%%%%%%%%%%%%%%%%%%%%%%%%%%%%%%%%%%%%%%%%%%%%%%%
Recall the definition and inversion formula for Dedekind's eta
function \cite[Theorem 3.1]{TMA}:
\begin{align*}
\eta(\tau)=e^{\pi i\tau/12}\prod_{m=1}^{\infty}(1-e^{2\pi i m\tau}),
 \qquad
\eta(-1/\tau)=\sqrt{-i\tau}\,\eta(\tau)\quad(\tau \in \mathbb{C}
\text{ and }\mathrm{Im}(\tau)>0).
\end{align*}
Let $\sigma_1(k)$ denote the sum of the positive divisors of $k$, and, for $t>0$, define
\begin{align*}
D(t)=\sum_{\substack{k\ge1\\k\ \rm{odd}}}\sigma_1(k)e^{-kt},\qquad W(t)=\log\mathcal{B}(e^{-t}).  
\end{align*}
Taking logarithms of both sides of \eqref{eq7.2}, we obtain
\begin{align}\label{eq7.77}
\log\mathcal{B}(q)&=\sum_{\substack{s\ge1\\s\ \rm{odd}}}\left(\log(1+q^s)-\log(1-q^s)\right)\nonumber\\
&=2\sum_{\substack{j\ge1\\j\ \rm{odd}}}
\sum_{\substack{r\ge1\\r\ \rm{odd}}}\frac{q^{jr}}r=2\sum_{\substack{k\ge1\\k\ \rm{odd}}}\frac{\sigma_1(k)}k q^k.
\end{align}
The last equality follows by taking $k=jr$: Both $j$ and $r$ are
odd precisely when $k$ is odd, and the coefficient of $q^k$ is $2\sum_{j\mid k}\frac1{k/j}=2\sigma_1(k)/k$. Substituting $q=e^{-t}$ into \eqref{eq7.77} and then differentiating, we obtain
%\begin{align}\label{eq7.8}
%W(t):=\log\mathcal{B}(e^{-t})=2\sum_{\substack{k\ge1\\k\ \rm{odd}}}\frac{\sigma_1(k)}ke^{-kt}.
%\end{align}
%Differentiating gives
\begin{align}\label{eq7.9}
W'(t)=-2D(t),\qquad W'''(t)=-2D''(t).
\end{align}
%%%%%%%%%%%%%%%%%%%%%%%%%%%%%%%%%
\begin{lemma}\label{lemma5.3a}
For $0<t\le1/100$,
\begin{align*}
D(t)\ge\frac{\pi^2}{16t^2}-1,\qquad D''(t)\le\frac{3\pi^2}{8t^4}+1.
\end{align*}
\end{lemma}
\begin{proof}
For $\tau=ijt/(2\pi)$ ($j\geq1$), substituting $\eta(\tau)$ and $\eta(-1/\tau)$ from the definition into the inversion formula yields
\begin{align}\label{eq7.6}
f_j(e^{-t})
 =\sqrt{\frac{2\pi}{jt}}\,
   \exp\left(\frac{jt}{24}-\frac{\pi^2}{6jt}\right)
   f_1(e^{-4\pi^2/(jt)}).
\end{align}
Applying \eqref{eq7.6} with $j=1,2,4$ to the factors in \eqref{eq7.2}, we obtain
\begin{align}\label{eq7.7}
\mathcal{B}(e^{-t})=\frac{f_2(e^{-t})^3}{f_1(e^{-t})^2f_4(e^{-t})}=\frac1{\sqrt2}e^{\pi^2/(8t)}\frac{f_2(e^{-\pi^2/t})^3}{f_1(e^{-\pi^2/t})f_4(e^{-\pi^2/t})^2}.   
\end{align}
Define 
\begin{align*}
\Omega(t)=\log\left(\frac{f_2(e^{-\pi^2/t})^3}{f_1(e^{-\pi^2/t})f_4(e^{-\pi^2/t})^2}\right).  
\end{align*}
Then, taking logarithms of both sides of \eqref{eq7.7}, we get
\begin{align}\label{eq7.10}
W(t)=\frac{\pi^2}{8t}-\frac12\log2+\Omega(t).
\end{align}
Combining \eqref{eq7.9} and \eqref{eq7.10}, we obtain
\begin{align}\label{eq7.88}
D(t)=\frac{\pi^2}{16t^2}-\frac12\Omega'(t),
 \qquad
D''(t)=\frac{3\pi^2}{8t^4}-\frac12\Omega'''(t).
\end{align}
It remains to bound $\Omega'(t)$ and $\Omega'''(t)$. Since
\begin{align*}
\log f_j(e^{-\pi^2/t})=\sum_{d\ge1}\log(1-e^{-\pi^2jd/t})=-\sum_{d\ge1}\sum_{r\ge1}\frac{e^{-\pi^2jdr/t}}r
  =-\sum_{m=1}^{\infty}\frac{\sigma_1(m)}m e^{-\pi^2jm/t},
\end{align*}
we have
\begin{align}\label{eq7O}
\Omega(t)=\sum_{m=1}^{\infty}w_me^{-\pi^2m/t},   
\end{align}
where
\begin{align}\label{eqw_m}
w_m=\frac{\sigma_1(m)}m
 -\mathbf1_{2\mid m}\cdot3\frac{\sigma_1(m/2)}{m/2}
 +\mathbf1_{4\mid m}\cdot 2\frac{\sigma_1(m/4)}{m/4}.    
\end{align}

We take $\sigma_1(x)=0$ if $x$ is not an integer. Set $u=\pi^2/t$. For $\gamma=2,3,4$, define
\begin{align*}
M_\gamma(u)=\sum_{m=1}^{\infty}m^{\gamma}e^{-mu}. 
\end{align*}
For $0<x\le1/2$, repeated application of the operator
$x\,d/dx$ to $\sum_{m=1}^{\infty}x^m=x/(1-x)$ gives
\begin{align}\label{eq7.11}
 \sum_{m=1}^{\infty}m^4x^m
 =\frac{x(1+11x+11x^2+x^3)}{(1-x)^5}\le x\,
 \frac{1+11/2+11/4+1/8}{(1-1/2)^5}
 =300x.    
\end{align}
Since $m^2\le m^4$ and $m^3\le m^4$ for $m\ge1$, taking
$x=e^{-u}<1/2$ in \eqref{eq7.11} implies
\begin{align}\label{eq7.12}
 M_\gamma(u)\le300e^{-u}
 \qquad(\gamma=2,3,4).
\end{align}
For every positive integer $r$, $\sigma_1(r)\le r^2$. It follows from \eqref{eqw_m} that $|w_m|\leq 3m$. For $t\leq1/100$, $u>900$. For the first derivative, \eqref{eq7O} gives
\begin{align*}
|\Omega'(t)|\le\frac{\pi^2}{t^2}\sum_{m=1}^{\infty}m|w_m|e^{-mu}\le\frac{3\pi^2}{t^2}M_2(u)=\frac{3u^2}{\pi^2}M_2(u)\le\frac{900}{\pi^2}u^2e^{-u}<100u^2e^{-u},   
\end{align*}
where the penultimate inequality follows from \eqref{eq7.12}. Similarly, for the third derivative, the triangle inequality yields
\begin{align*}
|\Omega'''(t)|
 &\le \frac{3\pi^6}{t^6}M_4(u)+\frac{18\pi^4}{t^5}M_3(u)+\frac{18\pi^2}{t^4}M_2(u)=\frac3{\pi^6}\left(u^6M_4(u)+6u^5M_3(u)+6u^4M_2(u)\right) \\
 &\le\frac{900}{\pi^6}(u^6+6u^5+6u^4)e^{-u}<\frac{100}{81}u^6\left(1+\frac6u+\frac6{u^2}\right)e^{-u}\leq 2u^6e^{-u},
\end{align*}
where the last inequality uses $u>900$.
For $u>0$, we have $e^u>\frac{u^4}{4!}$ and $e^u>\frac{u^8}{8!}$.
It follows that
\begin{align*}
|\Omega'(t)|&<100u^2e^{-u}<\frac{100\cdot4!}{u^2}=\frac{2400}{u^2}\le\frac{2400}{900^2}<1,\\
|\Omega'''(t)|&<2u^6e^{-u}<\frac{2\cdot8!}{u^2}=\frac{80640}{u^2}\le\frac{80640}{900^2}<1.
\end{align*}
Using these estimates in \eqref{eq7.88} completes the proof.
\end{proof}
%%%%%%%%%%%%%%%%%%%%%%_Lemma_5.3_%%%%%%%%%%%%%%%%%%%%%%%%%%%%%%%%%%%%%%%%%%%%%%%%
\begin{lemma}\label{lemma5.3}
For every integer $n\ge0$,
 \begin{align*}
 \beta(n)>\frac1{256}e^{(13/6)\sqrt n}.
 \end{align*}
\end{lemma}
\begin{proof}
Differentiating \eqref{eq7.77} gives the following recurrence
\begin{equation}\label{eq7.13}
n \beta(n)=2\sum_{\substack{1\le k\le n\\k\ \rm{odd}}}\sigma_1(k)\beta(n-k)\qquad(n\ge1).
\end{equation}
Set
\begin{align*}
\phi(n)=\frac1{256}e^{(13/6)\sqrt n}.
\end{align*}
We prove $\beta(n)>\phi(n)$ for all $n\geq0$ by induction. We verify using Mathematica that $\beta(n)>\phi(n)$ for $0\le n<20000$.
Suppose now that $n\ge20000$ and that $\beta(j)>\phi(j)$ for every $j<n$.
For $0\le x\le1$, putting $s=\sqrt{1-x}$ gives
\begin{align}\label{eq7.14}
 \sqrt{1-x}-1+\frac x2+\frac{x^2}{2}=\frac12s(1-s)^2(2+s)\ge0.
\end{align}
Applying \eqref{eq7.14} with $x=k/n$, and using $e^{-y}\ge1-y$, yields
\begin{align}\label{eq7.15}
\frac{\phi(n-k)}{\phi(n)}\ge e^{-13k/(12\sqrt{n})}\left(1-\frac{13k^2}{12n^{3/2}}\right)
\qquad(0\le k\le n).    
\end{align}
For $k>n$, the expression on the right is negative, because
$13\sqrt n/12>1$. Thus the lower-bound sum obtained from
\eqref{eq7.15} may be extended to all positive odd $k$.
For $n\geq20000$, $0<13/(12\sqrt{n})<1/100$, and therefore, \eqref{eq7.15} and Lemma \ref{lemma5.3a} (with $t=13/(12\sqrt{n})$) give
\begin{align}
\frac2n\sum_{\substack{1\le k\le n\\k\ \rm{odd}}}\sigma_1(k)\frac{\phi(n-k)}{\phi(n)}
&\ge\frac2nD\left(\frac{13}{12\sqrt{n}}\right)-\frac{13}{6n^{5/2}}D''\left(\frac{13}{12\sqrt{n}}\right)\nonumber\\
&\ge\frac{18\pi^2}{13^2}-\frac{6^4\pi^2}{13^3\sqrt n}-\frac2n-\frac{13}{6n^{5/2}}
>1,
\label{eq7.16}
\end{align}
as $13<6n^{3/2}$ for $n\ge20000$.
Using the induction hypothesis in \eqref{eq7.13} and employing \eqref{eq7.16}, we conclude that
\begin{align*}
\beta(n)>\frac2n\sum_{\substack{1\le k\le n\\k\ \rm{odd}}}
 \sigma_1(k)\phi(n-k)>\phi(n).
\end{align*}
This completes the proof.
\end{proof}
We are now ready to prove Proposition \ref{prop5.1}.
\begin{proof}[Proof of Proposition \ref{prop5.1}]
The ratio of the lower bound for $\beta(n)$ to the upper bound for $\alpha(n)$ obtained in Lemmas \ref{lemma5.2} and \ref{lemma5.3}, respectively, is
\begin{align*}
\frac{e^{(13/6)\sqrt n}/256}{2e^{(53/25)\sqrt n}}
 =\frac1{512}e^{(7/150)\sqrt n}>1,   
\end{align*}
for $n\ge20000$. Consequently $\beta(n)>\alpha(n)$ for all $n\ge20000$. For $28\le n\le 20000$, exact computation by Mathematica gives the same inequality. 
Finally, \eqref{eq7p} implies that for $n>27$,
\begin{align*}
\overline{\mathcal{P}}_e(n)>\overline{\mathcal{M}}_e(n).
\end{align*}
This completes the proof.
\end{proof}
%%%%%%%%%%%%%%%%%%%%%%%%%%%%%%%%%%%%%%%%%%%%%%%%%%%%%%%%%%%%%%%%%%%%%%%%%%%%%%%%%%%%%%%%%%%%%%%%%%%%%%%%%%%%%%%%%%%%%%%%%%%%%%%%%%%%%%%%%%%%%%%%%%%%%%%%%%%%%%%%%%%%%%%%%%%%%%%%%%%%%%%%%%%%%%%%%%%%%%%%%%%%%%%%%%%%%%%%%%%%%%%%%%%
\section{Concluding Remarks}
Let $\mathcal{M}_e(k,n)$ (resp., $\mathcal{M}_o(k,n)$) denote the number of partitions of $n$ with an even (resp., odd) number of distinct positive integers not exceeding the largest part and appearing less than $k$ times. Let $\mathcal{P}_e(k,n)$ (resp., $\mathcal{P}_o(k,n)$) denote the number of partitions of $n$ with an even (resp., odd) number of parts different from $k$. For $k=1$, Bhoria et al. conjectured \cite[Conjecture 2.7]{BES_2026} and the authors \cite[Corollary 1.5]{DS_2026} proved that $\mathcal{M}_e(1,n)>\mathcal{M}_o(1,n)$ for $n>34$. (In the previous sections, we suppress the first argument when $k=1$ and write $\mathcal{X}_\star(n)$ in place of $\mathcal{X}_\star(1,n)$, for $\mathcal{X}=\mathcal{M}, \mathcal{P}$ and $\star=e,o$.) In this article (Theorems \ref{Thm_main_2} and \ref{Thm_main_3}), we prove that
\begin{align*}
\mathcal{M}_e(1,n)&>\mathcal{P}_e(1,n)>\mathcal{P}_o(1,n)>\mathcal{M}_o(1,n)\quad(\text{for even } n>34),\\
\mathcal{M}_e(1,n)&>\mathcal{P}_o(1,n)>\mathcal{P}_e(1,n)>\mathcal{M}_o(1,n)\quad(\text{for odd } n>34).
\end{align*}
A natural question that arises is: What would the inequalities be for $k\geq2$?
\par A similar story is in the case of overpartitions. Let $\overline{\mathcal{M}}_e(k,n)$ (resp., $\overline{\mathcal{M}}_o(k,n)$) denote the number of overpartitions of $n$ with an even (resp., odd) number of distinct positive integers not exceeding the largest part and appearing less than $k$ times. Let $\overline{\mathcal{P}}_e(k,n)$ (resp., $\overline{\mathcal{P}}_o(k,n)$) denote the number of overpartitions $\overline{\lambda}$ of $n$ with $\psi_k(\overline{\lambda})$ even (resp., odd). Here, $\psi_k(\overline{\lambda})$ is the number of nonoverlined parts of $\overline{\lambda}$ with sizes in $\mathcal{S}_k$, plus its number of overlined parts when $k$ is even. For $k=1$, Bhoria et al. conjectured \cite[Conjecture 3.6]{BES_2026} and the authors \cite[Corollary 1.5]{DS_2026} proved that $\overline{\mathcal{M}}_e(1,n)>\overline{\mathcal{M}}_o(1,n)$ for $n>27$. In this article, we prove that for all $n>27$ (Theorem \ref{Thm_main_5}),
\begin{align*}
\overline{\mathcal{P}}_e(1,n)>\overline{\mathcal{M}}_e(1,n)>\overline{\mathcal{M}}_o(1,n)>\overline{\mathcal{P}}_o(1,n).
\end{align*}
The same question arises here: What would the inequalities be for $k\geq2$?

\bigskip
%%%%%%%%%%%%%%%%%%%%%%%%%%%%%%%%%%%%%%%
\noindent \textbf{Acknowledgements.}
This work was supported by
the National Natural Science Foundation of
China (12201177), the Natural Science Foundation of Hebei Province
(A2024205012, A2023205045), Hebei Key Talent Program (Overseas Returnee Platform) (A2026023), and
the 111 Center (D26018).
% The second author was supported by the Natural Science Foundation of Hebei Province (A2023205045) and the 111 Center (D26018).

\bigskip
%%%%%%%%%%%%%%%%%%%%%%%%%%%%%%%%%%%%%%%%%%%%%%%%%%%%%%%%%%%%%%%%%%%%%%%%%%%%%%%%%%%%%%%%%%%%%%%%%%%%%%%%%%%%%%%%%%%%%%%%%%%%%%%%%%%%%%%%%%%%%%%%%%%%%%%%%%%%%%%%%%%%%%%%%%%%%%%%%%%%%%%%%%%%%%%%%%%%%%%%%%%%%%%%%%%%%%%%%%%%%%%%%%%%%%%%%%%%%%%%%%%%%%%%%%%%%%%%%%%%%%%%%%%%%%%%%%%%%%%%%%%%%%%%%%%%%%%%%%%%%%%%%%%%%%%%%%%%%%%%%%%%%%%%%%%%%%%%%%%%%%%%%%%%%%%%%%%%%%%%%%%%%%%%%%%%%%%%%%%%%%%%%%%%%%%%%%%%%%%%%%%%%%%%%%%%%%%%%%%%%%%%%%%%%%%%%%%%%%%%%%%%%%%%%%%%%%%%%%%%%%%%%%%%%%%%%%%%%%%%%%%%%%%%%%%%%%%%%%%%%%%%%%%%%%%%%%%%%%%%%%%%%%%%%%%%%%%%%%%%%%%%%%%%%%%%%%%%%%%%%%%%%%%%%%%%%%%%%%%%%%%%%%%%%%%%%%%%%%%%%%%%%%%%%%%%%%%%%%%%%%%%%%%%%%%%%%%%%%%%%%%%%%%%%%%%%%%%%%%%%%%%%%%%%%%%%%%%%%%%%%%%%%%%%%%%%%%%%%%%%%%%%%%%%%%%%%%%%%%%%%%%%%%%%%%%%%%%%%%%%%%%%%%%%%%%%%%%%

%\noindent \textbf{Declaration of AI use.}
%ChatGPT (OpenAI) was used to assist in reducing the error term and lowering the lower bound of $n$ in the circle-method proof of Theorem \ref{Thm_main} in Section \ref{sect-missing}. All AI-assisted calculations, arguments, and suggestions were independently verified and revised by the authors, who take full responsibility for the content of the article.

%%%%%%%%%%%%%%%%%%%%%%%%%%%%%%%%%%%%%%%%%%%%%%%
%%%%%%%%%%%%%%%%%%%%%%%%%%%%%%%%%%%%%%%%%%%%%%%
%%%%%%%%%%%%%%%%%%%%%%%%%%%%%%%%%%%%%%%%%%%%%%%

\end{document}